\documentclass[11pt]{article}
\usepackage[ansinew]{inputenc}
\usepackage{array}
\usepackage{color}
\usepackage{amsmath}
\usepackage{amsxtra}
\usepackage{amstext}
\usepackage{amssymb}
\usepackage{latexsym}
\usepackage{amsfonts,cite}
\usepackage{graphics,epstopdf}
\usepackage{epsfig}

\begin{document}

\sloppy
\newcommand{\proof}{{\it Proof.~}}
\newtheorem{thm}{Theorem}[section]
\newtheorem{cor}{Corollary}[section]
\newtheorem{lem}{Lemma}[section]
\newtheorem{prop}[thm]{Proposition}
\newtheorem{eg}{Example}[section]
\newtheorem{defn}{Definition}[section]

\newtheorem{rem}{Remark}[section]
\numberwithin{equation}{section}

\thispagestyle{empty}
\parindent=0mm
\begin{center}
{\Large{\bf Generalization of Sz\'{a}sz operators involving\\
 multiple Sheffer polynomials}}\\
\vspace{0.3cm}
{\bf Mahvish Ali$^{1, \star}$ and Richard B. Paris$^{2}$}\\
\vspace{0.2cm}
$^{1}$Department of Applied Sciences and Humanities, Faculty of Engineering and Technology
Jamia Millia Islamia (A Central University), New Delhi-110025, India\\
\vspace{0.2cm}
$^{2}$Department of Mathematics, Abertay University, Dundee DD1 1HG, UK\\

\footnote{$^{\star}$Corresponding author.}
\footnote{Emails:}
\footnote{mahvishali37@gmail.com (Mahvish Ali)}
\footnote{r.paris@abertay.ac.uk (Richard B. Paris)}

\end{center}

\parindent=0mm

\vspace{0.5cm}
\noindent
{\bf Abstract:}~~The present work deals with the mathematical investigation of some generalizations of the Sz\'{a}sz operators. In this work, the multiple Sheffer polynomials are introduced. The generalization of Sz\'{a}sz operators involving multiple Sheffer polynomials are considered. Convergence properties of these operators are verified with the help of the universal Korovkin-type property and the order of approximation is calculated by using classical modulus of continuity. The theoretical results are exemplified choosing the special cases of multiple Sheffer polynomials.\\

\noindent
{\bf{\em Keywords:}}~Sz\'{a}sz operators, modulus of continuity, rate of convergence, multiple Sheffer polynomials.\\

\noindent
{\bf {\em 2010 MSC:}}~41A10, 41A25, 41A36, 33C45, 33E20.

\parindent=8mm

\section{Introduction}

The positive approximation processes discovered by Korovkin \cite{Korovkin} play a central role and arise in a natural way in many problems connected with functional analysis, harmonic analysis, measure theory, partial differential equations and probability theory. In 1953, P. P. Korovkin \cite{Korovkin} discovered perhaps, the  most powerful and at the same time, simplest criterion in order to decide whether a given sequence $(K_{n})_{n\in \mathbb{N}}$ of positive linear operators on the space $C[0,1]$ is an approximation process, i.e. $K_{n}(f)\rightarrow f$ uniformly on $[0,1]$ for every $f \in C[0,1]$. Starting with this result, a considerable number of mathematicians have extended Korovkin's theorem to other function spaces or, more generally, to abstract spaces, such as Banach lattices, Banach algebras, Banach spaces and so on. Korovkin's work, in fact, delineated a new theory that may be called Korovkin-type approximation theory.

One of the well-known examples of positive linear  operators is Szasz operators \cite{Szasz}.
Sz\'{a}sz \cite{Szasz} introduced the following positive  linear operators:
\begin{equation}\label{Szasz}
S_{n}(f;x):=e^{-nx}\sum\limits_{k=0}^{\infty}\frac{(nx)^{k}}{k!}f\left(\frac{k}{n}\right),
\end{equation}
where $x\ge 0$ and $f\in C[0,\infty)$ whenever the above sum converges. In recent years, there is an increasing interest to study linear positive operators based on certain polynomials, such as Appell polynomials, Sheffer polynomials, and Boas-Buck polynomials. Jakimovski and Leviatan \cite{JakiLevi} obtained a generalization of Szasz operators by means of the Appell polynomials defined as follows:
\begin{equation}
P_{n}(f;x):=\frac{e^{-nx}}{g(1)}\sum\limits_{k=0}^{\infty}p_{k}(nx)f\left(\frac{k}{n}\right),
\end{equation}
where $p_k(x), ~k\ge 0$ are the Appell polynomials defined by $g(t)e^{xt}=\sum\limits_{k=0}^{\infty}p_k(x)\frac{t^k}{k!}$ and $g(t)=\sum\limits_{k=0}^{\infty}a_k\frac{t^k}{k!}$ is an analytic function in the disk $|t|<R$, $R>1$ and $g(1)\neq 0$. For $g(t)=1$, we obtain the Sz\'{a}sz operators \eqref{Szasz}.

Ismail \cite{Ismail} presented a generalization of Szasz and Jakimovski and Leviatan operators by using the Sheffer polynomials $s_{k}(x)$ as
\begin{equation}\label{Ismail}
T_{n}(f;x):=\frac{e^{-nxH(1)}}{g(1)}\sum\limits_{k=0}^{\infty}s_{k}(nx)f\left(\frac{k}{n}\right), \quad n\in \mathbb{N},
\end{equation}
where $s_k(x), ~k\ge 0$ are the Sheffer polynomials defined by $g(t)e^{xH(t)}=\sum\limits_{k=0}^{\infty}s_k(x)\frac{t^k}{k!}$ and $g(t)$ are defined as above, $H(t)=\sum\limits_{k=0}^{\infty}h_k\frac{t^k}{k!}$ is an analytic function in the disk $|t|<R$, $R>1$ and $g(1)\neq 0$; $H'(1)=1$. For $g(t)=1$ and $H(t)=t,$ \eqref{Ismail} reduces to \eqref{Szasz}.

Recently, research has been undertaken in an attempt at generalization of the Sz\'{a}sz operators associated with  multiple polynomial sets. By a multiple polynomial system we mean a set of polynomials $\{p_{k_1,k_2}(x)\}_{k_1,k_2=0}^{\infty}$ with degree $k_1+k_2$, $k_1, k_2\ge 0$.

In \cite{Lee}, Lee defined the multiple Appell polynomials and found several equivalent conditions for this class of polynomials. A multiple polynomial set $\{p_{k_1,k_2}(x)\}_{k_1,k_2=0}^{\infty}$ is called multiple Appell if there exists a generating function of the form:
\begin{equation}\label{Multiple-Appell}
A(t_1,t_2)e^{x(t_1+t_2)}=\sum\limits_{k_1=0}^{\infty}\sum\limits_{k_2=0}^{\infty}p_{k_1,k_2}(x)\frac{t_{1}^{k_1}t_{2}^{k_2}}{k_{1}!~k_{2}!},
\end{equation}
where
$A(t_1,t_2)$ is defined as
\begin{equation}\label{A(t,t)A}
A(t_1,t_2)=\sum\limits_{k_1=0}^{\infty}\sum\limits_{k_2=0}^{\infty}a_{k_1,k_2}\frac{t_{1}^{k_1}t_{2}^{k_2}}{k_{1}!~k_{2}!}.
\end{equation}
The generalization of the Sz\'{a}sz operators involving multiple Appell polynomials has been studied in \cite{Varma}. The Jakimovski-Leviatan-Durrmeyer type operators involving multiple Appell polynomial are introduced in \cite{ansari} and investigate Korovkin type approximation theorem and rate of convergence. Some properties of the generalized Sz\'{a}sz operators by multiple Appell polynomials are given in \cite{Braha} taking into consideration the power summability method.

Inspired by the above works, we construct the generalization of the Sz\'{a}sz operators involving the multiple Sheffer polynomials. The paper is organized as follows: In Section 2, the multiple Sheffer polynomials are introduced and the positive linear operators $G_n(f;x)$ involving multiple Sheffer polynomials are constructed. In Section 3, some auxiliary results for the operators $G_n(f;x)$ are given. Section 4 is discusses approximation properties of the operators $G_n(f;x)$ and the convergence theorem. In order to show the relevance of the results, in the last section some numerical examples are given.

\section{Multiple Sheffer polynomials}

In this section, we introduce the multiple Sheffer polynomials as follows:

\begin{defn}
The multiple Sheffer polynomials set $\{S_{k_1,k_2}(x)\}_{k_1,k_2=0}^{\infty}$ possesses the following generating function:
\begin{equation}\label{Multiple-Sheffer}
A(t_1,t_2)e^{xH(t_1,t_2)}=\sum\limits_{k_1=0}^{\infty}\sum\limits_{k_2=0}^{\infty}S_{k_1,k_2}(x)\frac{t_{1}^{k_1}t_{2}^{k_2}}{k_{1}!~k_{2}!},
\end{equation}
where
$A(t_1,t_2)$ and $H(t_1,t_2)$ have series expansions of the form
\begin{equation}\label{A(t,t)}
A(t_1,t_2)=\sum\limits_{k_1=0}^{\infty}\sum\limits_{k_2=0}^{\infty}a_{k_1,k_2}\frac{t_{1}^{k_1}t_{2}^{k_2}}{k_{1}!~k_{2}!}
\end{equation}
and
\begin{equation}\label{H(t,t)}
H(t_1,t_2)=\sum\limits_{k_1=0}^{\infty}\sum\limits_{k_2=0}^{\infty}h_{k_1,k_2}\frac{t_{1}^{k_1}t_{2}^{k_2}}{k_{1}!~k_{2}!},
\end{equation}
respectively with the conditions that
\begin{equation}
A(0,0)=a_{0,0}\neq 0, \quad  H(0,0)=h_{0,0}\neq 0.
\end{equation}
\end{defn}

For $A(t_1,t_2)=e^{\frac{\delta}{2}(t_1+t_2)^2+\alpha_1 t_1+\alpha_2 t_2}$ and $H(t_1,t_2)=\delta(t_1+t_2)$, the multiple Sheffer polynomials $S_{k_1,k_2}(x)$ become the multiple Hermite polynomials $H_{k_1,k_2}^{(\alpha_1,\alpha_2)}(x)$ defined by the generating function \cite{Lee2}:
\begin{equation}\label{multiHermite}
e^{\frac{\delta}{2}(t_1+t_2)^2+\alpha_1 t_1+\alpha_2 t_2+\delta(t_1+t_2)x}=\sum\limits_{k_1=0}^{\infty}\sum\limits_{k_2=0}^{\infty}H_{k_1,k_2}^{(\alpha_1,\alpha_2)}(x)\frac{t_{1}^{k_1}t_{2}^{k_2}}{k_{1}!~k_{2}!},
\end{equation}
with $\delta<0$ and $\alpha_1 \neq \alpha_2$.
Note that if we take $t_2 = 0$, we get the generating function for the classical Hermite polynomials.\\

For $A(t_1,t_2)=\frac{1}{(1-t_1-t_2)^{\alpha+1}}$ and $H(t_1,t_2)=\frac{\beta_1 t_1+\beta_2 t_2}{1-t_1-t_2}$, the multiple Sheffer polynomials $S_{k_1,k_2}(x)$ become the multiple Laguerre polynomials $L_{k_1,k_2}^{(\alpha; \beta_1, \beta_2)}(x)$ defined by the generating function \cite{Lee2}:
\begin{equation}\label{multiLaguerre}
\frac{1}{(1-t_1-t_2)^{\alpha+1}}e^{\left(\frac{\beta_1 t_1+\beta_2 t_2}{1-t_1-t_2}\right)x}=\sum\limits_{k_1=0}^{\infty}\sum\limits_{k_2=0}^{\infty}L_{k_1,k_2}^{(\alpha; \beta_1, \beta_2)}(x)\frac{t_{1}^{k_1}t_{2}^{k_2}}{k_{1}!~k_{2}!},
\end{equation}
with $\alpha >-1$ and $\beta_1 \neq \beta_2$.
Note that if we take $t_2 = 0$, we get the generating function for the classical Laguerre polynomials.\\

Now, we construct the positive linear operators involving multiple Sheffer polynomials $S_{k_1,k_2}(x)$.
Throughout the paper, the following abbreviations for the partial derivatives will be used:
\begin{equation}
\frac{\partial A}{\partial t_i}=A_{t_i}, \quad \frac{\partial H}{\partial t_i}=H_{t_i}
\end{equation}
and
\begin{equation}
  \frac{\partial^2 A}{\partial t_i \partial t_j}=A_{t_i, t_j}, \quad \frac{\partial^2 H}{\partial t_i \partial t_j}=H_{t_i, t_j}, \quad i,j=1,2.
\end{equation}

Let us consider the following conditions on $S_{k_1,k_2}(x)$:
\begin{equation}\label{restrictions}
\begin{split}
(i)&~S_{k_1,k_2}(x)\ge 0, k_1,k_2\in \mathbb{N},\\
(ii)&~A(1,1)\neq 0, H_{t_1}(1,1)=1, H_{t_2}(1,1)=1,\\
(iii)&~\eqref{Multiple-Sheffer}, \eqref{A(t,t)}~{\text{ and}}~ \eqref{H(t,t)} ~{\text{converge~ for}}~ |t_1|<R_1, ~|t_2|<R_2 ~(R_1, R_2)>1.\\
\end{split}
\end{equation}

Under the above conditions, we construct a positive linear operator involving multiple Sheffer polynomials for $x\in [0,\infty)$ as follows:
\begin{equation}\label{MultipleShef-Operator}
G_{n}(f;x)=\frac{e^{-\frac{nx}{2}H(1,1)}}{A(1,1)}\sum\limits_{k_1=0}^{\infty}\sum\limits_{k_2=0}^{\infty}\frac{S_{k_1,k_2}(\frac{nx}{2})}{k_{1}!~k_{2}!}
f\left(\frac{k_{1}+k_2}{n}\right),
\end{equation}
provided that the right-hand side of \eqref{MultipleShef-Operator} exists.

\begin{rem}
It is to be noted that if we consider $H(t_1,t_2)=t_1+t_2$, the multiple Sheffer polynomials $S_{k_1,k_2}(x)$ become the multiple Appell polynomials $p_{k_1,k_2}(x)$ defined by \eqref{Multiple-Appell}. Therefore, for $H(t_1,t_2)=t_1+t_2$ the operators \eqref{MultipleShef-Operator} reduce to the operators involving multiple Appell polynomials defined by Varma \cite{Varma}.
\end{rem}

\begin{rem}
For $A(t_1,t_2)=1$, $H(t_1,t_2)=t_1+t_2$ and then $t_2=0$, the multiple Sheffer polynomials $S_{k_1,k_2}(x)$ reduce to $x^k$.  Therefore, the operators \eqref{MultipleShef-Operator} reduce to the Sz\'{a}sz operators \eqref{Szasz}.
\end{rem}

\section{Auxiliary results}

Note that throughout the paper we will assume that the operators  $G_{n}(f;x)$ are positive and we use the following test functions:
$$e_i(x) = x^i; \quad  i\in \{0,1,2\}.$$

\begin{lem}\label{1}
From  \eqref{Multiple-Sheffer}, the following hold:
\begin{eqnarray*}
\sum\limits_{k_1=0}^{\infty}\sum\limits_{k_2=0}^{\infty}\frac{S_{k_1,k_2}(\frac{nx}{2})}{k_{1}!~k_{2}!}&=&A(1,1)e^{\frac{nx}{2}H(1,1)},\label{e0}\\
\sum\limits_{k_1=0}^{\infty}\sum\limits_{k_2=0}^{\infty}(k_1+k_2)\frac{S_{k_1,k_2}(\frac{nx}{2})}{k_{1}!~k_{2}!}&=&\left[nxA(1,1)+A_{t_1}(1,1)+A_{t_2}(1,1)\right]e^{\frac{nx}{2}H(1,1)},\label{e1}\\
\sum\limits_{k_1=0}^{\infty}\sum\limits_{k_2=0}^{\infty}(k_1+k_2)^2\frac{S_{k_1,k_2}(\frac{nx}{2})}{k_{1}!~k_{2}!}&
=&\left[n^2x^2A(1,1)+nx\left(\frac{A(1,1)H_{t_1,t_1}(1,1)+A(1,1)H_{t_2,t_2}(1,1)}{2}\right.\right.\\
&&\left.+A(1,1)H_{t_1,t_2}(1,1)+2A_{t_1}(1,1)+2A_{t_2}(1,1)+A(1,1)\vphantom{\frac{A(1,1)H_{t_1,t_1}(1,1)+A(1,1)H_{t_2,t_2}(1,1)}{2}}\right)\\
&&+A_{t_1,t_1}(1,1)+A_{t_2,t_2}(1,1)+2A_{t_1,t_2}(1,1)\\
&&\left.+A_{t_1}(1,1)+A_{t_2}(1,1)\vphantom{\frac{A(1,1)H_{t_1,t_1}(1,1)+A(1,1)H_{t_2,t_2}(1,1)}{2}}\right]e^{\frac{nx}{2}H(1,1)}.\label{e2}\\
\end{eqnarray*}
\end{lem}

\noindent
\begin{proof}
Taking the partial derivatives of the generating function \eqref{Multiple-Sheffer} with respect to $t_1$ and $t_2$ and then letting $t_1=1,~ t_2=1$, we obtain Lemma \ref{1}.
\end{proof}

\begin{lem}\label{2}
The operators $G_n$ defined by \eqref{MultipleShef-Operator} satisfy
\begin{eqnarray}
G_n(e_0;x)&=&1\label{E00}\\
G_n(e_1;x)&=&x+\frac{A_{t_1}(1,1)+A_{t_2}(1,1)}{n~A(1,1)}\label{E11}\\
G_n(e_2;x)&=&x^2+\frac{x}{n}\left(1+\frac{H_{t_1,t_1}(1,1)+H_{t_2,t_2}(1,1)}{2}+\frac{2A_{t_1}(1,1)+2A_{t_2}(1,1)}{A(1,1)}+H_{t_1,t_2}(1,1)\right)\nonumber\\
&&+\frac{A_{t_1,t_1}(1,1)+A_{t_2,t_2}(1,1)+2A_{t_1,t_2}(1,1)+A_{t_1}(1,1)+A_{t_2}(1,1)}{n^2 ~A(1,1)}.\label{E22}
\end{eqnarray}
\end{lem}

\noindent
\begin{proof}
In view of operators \eqref{MultipleShef-Operator} and Lemma \ref{1}, the proof of Lemma \ref{2} is
straightforward.\\
\end{proof}

By using the linearity property of the operators \eqref{MultipleShef-Operator} and Lemma \ref{2}, it follows that:
\begin{eqnarray}
G_n((e_1-x);x)&=&\frac{A_{t_1}(1,1)+A_{t_2}(1,1)}{n~A(1,1)}\\
G_n((e_1-x)^2;x)&=&\frac{x}{n}\left(1+\frac{H_{t_1,t_1}(1,1)+H_{t_2,t_2}(1,1)}{2}+H_{t_1,t_2}(1,1)\right)\nonumber\\
&&+\frac{A_{t_1,t_1}(1,1)+A_{t_2,t_2}(1,1)+2A_{t_1,t_2}(1,1)+A_{t_1}(1,1)+A_{t_2}(1,1)}{n^2 ~A(1,1)}.
\end{eqnarray}

\section{Approximation results}

In this section, we state our main theorem with the help of the universal Korovkin-type property and calculate the order of approximation by modulus of continuity.

P. P. Korovkin \cite{Korovkin} has proved some remarkable results concerning the convergence of sequences $(K_{n}(f,x))_{n=1}^{\infty}$, where $K_{n}(f,x)$ are positive linear operators. For example, if $K_{n}(f,x)$ converges uniformly to $f$ in the particular cases $f(t)\equiv 1$, $f(t) \equiv t$, $f(t) \equiv t^2$, then it does so for every continuous real $f$.

First, we recall the following definitions and lemmas:\\

\begin{defn}
Let $f\in \hat{C}[0,\infty)$ and $\delta > 0$. The modulus of continuity $w(f;\delta)$ of the function f is defined by
\begin{equation}
w(f;\delta):=  \sup_{\substack{x,y\in[0,\infty)\\|x-y|\leq \delta}}|f(x)-f(y)|,
\end{equation}
where $\hat{C} [0,\infty)$ is the space of uniformly continuous functions on $[0,\infty)$. Then, for any $\delta>0$ and each $x\in [0,\infty)$, it is well known that one can write
\begin{equation}\label{modulus}
|f(x)-f(y)|\le w(f;\delta)\left(\frac{|x-y|}{\delta}+1\right).
\end{equation}

If $f$ is uniformly continuous on $[0,\infty)$, then it is necessary and sufficient that
$$\lim\limits_{\delta\rightarrow 0} w(f, \delta)= 0.$$
\end{defn}

\begin{defn}
The second modulus of continuity of the function $f\in C_{B}[0,\infty)$ is defined by
\begin{equation}
w_{2}(f;\delta):= \sup_{\substack{0<t\le \delta}}||f(.+2t)-2f(.+t)+f(.)||_{C_{B}},
\end{equation}
where $C_{B}[0,\infty)$ is the class of real-valued functions defined on $[0,\infty)$, which are bounded and uniformly continuous with the norm
\begin{equation}
||f||_{C_{B}}=\sup_{x\in [0,\infty)}|f(x)|.
\end{equation}
\end{defn}

Let us define the class $E$ as follows:
$$E:=\left\{f: x\in [0,\infty), ~\frac{f(x)}{1+x^2}~{\text{is~convergent~as}}~n\rightarrow \infty\right\}$$

\begin{lem}\label{Gav}
(Gavrea and Ra\c{s}a \cite{Gavrea}). Let $ g\in C^{2}[0,a]$ and $(K_{n})_{n\geq 0}$ be a sequence of positive linear operators with the property $K_{n}(1;x)=1$. Then,
\begin{equation}
|K_{n}(g;x)-g(x)|\leq||g'|| \sqrt{K_{n}((s-x)^{2};x)}+\frac{1}{2}||g''||~ K_{n}((s-x)^{2};x).
\end{equation}
\end{lem}

\begin{lem}\label{Zhuk}
(Zhuk \cite{Zhuk}). Let $f\in C[a,b]~~ \text{and}~~ h \in(0,\frac{a-b}{2})$. Let $f_{h}$ be the second-order Steklov function attached to the function f. Then, the following inequalities are satisfied:
\begin{equation}
(i)~ ||f_{h}-f||\leq\frac{3}{4}w_{2}(f;h),\label{2.10}
\end{equation}
\begin{equation}
(ii)~ ||f''_{h}||\leq \frac{3}{2h^{2}}w_{2}(f;h).\label{2.11}
\end{equation}
\end{lem}

\begin{thm}
Let $f\in C[0,\infty)\cap E$. Then
\begin{equation}
\lim_{n\rightarrow \infty}G_{n}(f;x)=f(x)
\end{equation}
uniformly on each compact subset of $[0,\infty)$.
\end{thm}

\noindent
\begin{proof}
In view of Lemma \ref{2}, it follows that
\begin{equation}\label{K}
\lim_{n\rightarrow \infty}G_{n}(e_i,x)=x^i, \quad i=0, 1, 2
\end{equation}
uniformly on each compact subset of $[0,\infty)$. Application of Korovkin's theorem to \eqref{K} then establishes the desired result.\\
\end{proof}

Next, we obtain the order of approximation of the operators $G_n(f;x)$.

\begin{thm}
Let $f\in \hat{C}_{E}[0,\infty)$. Then the operators $G_n(f;x)$ satisfy the following inequality:
\begin{equation}\label{RC}
|{G_n(f;x)-f(x)}|\leqslant \{1+\sqrt{\lambda_n(x)}\}w(f;\sqrt{n}),
\end{equation}
where
$$\lambda_n(x)=x\left(1+\frac{H_{t_1,t_1}(1,1)+H_{t_2,t_2}(1,1)}{2}+H_{t_1,t_2}(1,1)\right)$$
$$+\frac{A_{t_1,t_1}(1,1)+A_{t_2,t_2}(1,1)+2A_{t_1,t_2}(1,1)+A_{t_1}(1,1)+A_{t_2}(1,1)}{n ~A(1,1)}.$$
\end{thm}

\noindent
\begin{proof}
In view of the modulus of continuity \eqref{modulus}, we have
\begin{eqnarray}
|{G_n(f;x)-f(x)}|&\leqslant & \frac{e^{-\frac{nx}{2}H(1,1)}}{A(1,1)}\sum\limits_{k_1=0}^{\infty}\sum\limits_{k_2=0}^{\infty}\frac{S_{k_1,k_2}(\frac{nx}{2})}{k_{1}!~k_{2}!}
\left|f\left(\frac{k_{1}+k_2}{n}\right)-f(x)\right|\nonumber\\
&\leqslant &w(f;\delta)\left\{ 1+ \frac{1}{\delta}\frac{e^{-\frac{nx}{2}H(1,1)}}{A(1,1)} \sum\limits_{k_1=0}^{\infty}\sum\limits_{k_2=0}^{\infty}\frac{S_{k_1,k_2}(\frac{nx}{2})}{k_{1}!~k_{2}!}\left|\frac{k_{1}+k_2}{n}-x\right|\right\}\label{3}
\end{eqnarray}

By considering the Cauchy-Schwarz inequality, we find
\begin{eqnarray}
&&\sum\limits_{k_1=0}^{\infty}\sum\limits_{k_2=0}^{\infty}\frac{S_{k_1,k_2}(\frac{nx}{2})}{k_{1}!~k_{2}!}\left|\frac{k_{1}+k_2}{n}-x\right|\nonumber\\
&\leqslant &\sqrt{A(1,1)e^{xH(1,1)}}
\left\{\sum\limits_{k_1=0}^{\infty}\sum\limits_{k_2=0}^{\infty}\frac{S_{k_1,k_2}(\frac{nx}{2})}{k_{1}!~k_{2}!}\left(\frac{k_{1}+k_2}{n}-x\right)^2 \right\}^{\frac{1}{2}}\nonumber\\
&&= A(1,1)e^{xH(1,1)}\left\{\frac{x}{n}\left(1+\frac{H_{t_1,t_1}(1,1)+H_{t_2,t_2}(1,1)}{2}+H_{t_1,t_2}(1,1)\right)\right.\nonumber\\
&&+\left.\frac{A_{t_1,t_1}(1,1)+A_{t_2,t_2}(1,1)+2A_{t_1,t_2}(1,1)+A_{t_1}(1,1)+A_{t_2}(1,1)}{n^2 ~A(1,1)} \right\}^{\frac{1}{2}}.\label{CS}
\end{eqnarray}

Use of \eqref{CS} in \eqref{3} then leads to
\begin{eqnarray}
|{G_n(f;x)-f(x)}|&\leqslant & w(f;\delta)\left\{ 1+ \frac{1}{\delta}\left[\frac{x}{n}\left(1+\frac{H_{t_1,t_1}(1,1)+H_{t_2,t_2}(1,1)}{2}+H_{t_1,t_2}(1,1)\right)\right.\right.\nonumber\\
&&\left.\left.+\frac{A_{t_1,t_1}(1,1)+A_{t_2,t_2}(1,1)+2A_{t_1,t_2}(1,1)+A_{t_1}(1,1)+A_{t_2}(1,1)}{n^2 ~A(1,1)}\right]^{\frac{1}{2}} \right\}.\nonumber\\
\end{eqnarray}

On choosing $\delta=\sqrt{n}$, assertion \eqref{RC} follows.\\
\end{proof}

\begin{thm}
For $f\in C[0,\alpha]$, the following inequality is satisfied:
\begin{equation}
|G_{n}(f;x)-f(x)| \leq \frac{2}{\alpha}||f||~l^{2}+\frac{3}{4}(\alpha+2+l^{2})w_{2}(f;l),\label{2.29}
\end{equation}
where
\begin{equation}
 l:=l_{n}(x)=\sqrt[4]{G_{n}((e_1-x)^{2};x)}\nonumber
\end{equation}
and the second-order modulus of continuity is given by $w_{2}(f;l)$ with the norm $||f||=\max_{\substack{x\in [a,b]}}|f(x)|$
\end{thm}

\noindent
\begin{proof}
 Let $f_{l}$ be the second-order Steklov function attached to the function $f$. In view of the identity \eqref{E00}, we have
\begin{eqnarray}
|G_{n}(f;x)-f(x)|&\leq &|G_{n}(f-f_{l};x)|+ |G_{n}(f_{l};x)-f_{l}(x)|+|f_{l}(x)-f(x)|\nonumber\\
&\leq& 2||f_{l}-f|| + |G_{n}(f_{l};x)-f_{l}(x)|,
\end{eqnarray}
which on using inequality \eqref{2.10} becomes
\begin{equation}
|G_{n}(f;x)-f(x)|\leq\frac{3}{2}w_{2}(f;l)+ |G_{n}(f_{l};x)-f_{l}(x)|.\label{2.31}
\end{equation}

Taking into account that $f_{l} \in C^{2}[0,\alpha]$, from Lemma \ref{Gav}, it follows that
\begin{equation}
|G_{n}(f_{l};x)-f_{l}(x)|\leq ||f'_{l}||  \sqrt{G_{n}((s-x)^{2};x)}+\frac{1}{2}~||f''_{l}||~ G_{n}((s-x)^{2};x),
\end{equation}
which in view of inequality \eqref{2.11} becomes
\begin{equation}
|G_{n}(f_{l};x)-f_{l}(x)|\leq ||f'_{l}||  \sqrt{G_{n}((s-x)^{2};x)}+\frac{3}{4l^2}w_{2}(f;l)G_{n}((s-x)^{2};x).\label{2.33}
\end{equation}

Further, the Landau inequality
\begin{equation}
||f' _{l}|| \leq \frac{2}{\alpha} ||f_{l}|| +\frac{ \alpha}{2} ||f'' _{l}||,\nonumber
\end{equation}
and Lemma \ref{Zhuk} gives
\begin{equation}
||f' _{l}||\leq \frac{2}{\alpha}||f||+\frac{3 \alpha}{4l^{2}}w_{2}(f;l).\label{2.34}
\end{equation}

Using inequality  \eqref{2.34} in inequality \eqref{2.33} and taking $ l=\sqrt[4]{{G_{n}((s-x)^{2});x}}$, we find
\begin{equation}
|G_{n}(f_{l};x)-f_{l}(x)|\leq \frac{2}{\alpha}||f||~ l^{2}+ \frac{3}{4}(\alpha+l^2)w_{2}(f;l).\label{2.35}
\end{equation}
Use of the inequality \eqref{2.35} in \eqref{2.31}, then establishes the assertion \eqref{2.29}.\\
\end{proof}

\section{Examples}

\begin{eg}
{\bf The case $\mathbf{\{(x+1)^{k_1+k_2}\}_{k_1,k_2=0}^{\infty}}$}
\end{eg}

By taking $A(t_1,t_2)=e^{t_1+t_2}$ and $H(t_1,t_2)=t_1+t_2$, the multiple Sheffer polynomials $S_{k_1,k_2}(x)$ reduce to the polynomials $\{(x+1)^{k_1+k_2}\}_{k_1,k_2=0}^{\infty}$ defined by the generating function:
\begin{equation}\label{xk}
e^{(x+1)(t_1+t_2)}=\sum\limits_{k_1=0}^{\infty}\sum\limits_{k_2=0}^{\infty}(x+1)^{k_1+k_2}\frac{t_{1}^{k_1}t_{2}^{k_2}}{k_{1}!~k_{2}!}.
\end{equation}
It is clear that $(x+1)^{k_1+k_2}\ge 0$ for $x\in (0,\infty]$, $A(1,1)=e^2\neq 0$ and $H_{t_1}(1,1)=H_{t_2}(1,1)=1$. Therefore, the conditions \eqref{restrictions} are satisfied.\\

The generalized Sz\'{a}sz operators \eqref{MultipleShef-Operator} involving the function $\{(x+1)^{k_1+k_2}\}_{k_1,k_2=0}^{\infty}$ are obtained as follows:
\begin{equation}\label{Gxk}
G_{n}(f;x)=e^{-2-nx}
\sum\limits_{k_1=0}^{\infty}\sum\limits_{k_2=0}^{\infty}
\frac{\left(\frac{nx}{2}+1\right)^{k_1+k_2}}{k_{1}!~k_{2}!}
f\left(\frac{k_{1}+k_2}{n}\right).
\end{equation}

\begin{lem}\label{G1L1}
For the operators $G_{n}(f;x)$ given by \eqref{Gxk} and $x\in [0,\infty)$, we have
\begin{eqnarray}
  G_{n}(e_0;x) &=& 1 \\
 G_{n}(e_1;x) &=& x+\frac{2}{n} \\
 G_{n}(e_2;x) &=& x^{2}+\frac{5x}{n}+\frac{6}{n^{2}}.
\end{eqnarray}
\end{lem}

\begin{lem}
For the operators $G_{n}(f;x)$ given by \eqref{Gxk} and $x\in [0,\infty)$, the following identities are satisfied:
\begin{eqnarray}
 G_{n}((e_1-1);x) &=& \frac{2}{n},\\
 G_{n}((e_1-x)^{2};x) &=&  \frac{x}{n}+\frac{6}{n^{2}}.\label{G1L2}
\end{eqnarray}
\end{lem}

For $n=20, 30, 50$, the convergence of the operators \eqref{Gxk} to the functions
\begin{subequations}
\begin{eqnarray}
f(x)&=&\left(x-\frac{1}{2}\right)\left(x-\frac{1}{3}\right)\label{f1}\\
f(x)&=&-4x^{3}\label{f2}
\end{eqnarray}
\end{subequations}
is illustrated in Figs. 1 and 2, respectively.\\

\newpage
\begin{figure}[htb]

\begin{center}
\epsfig{file=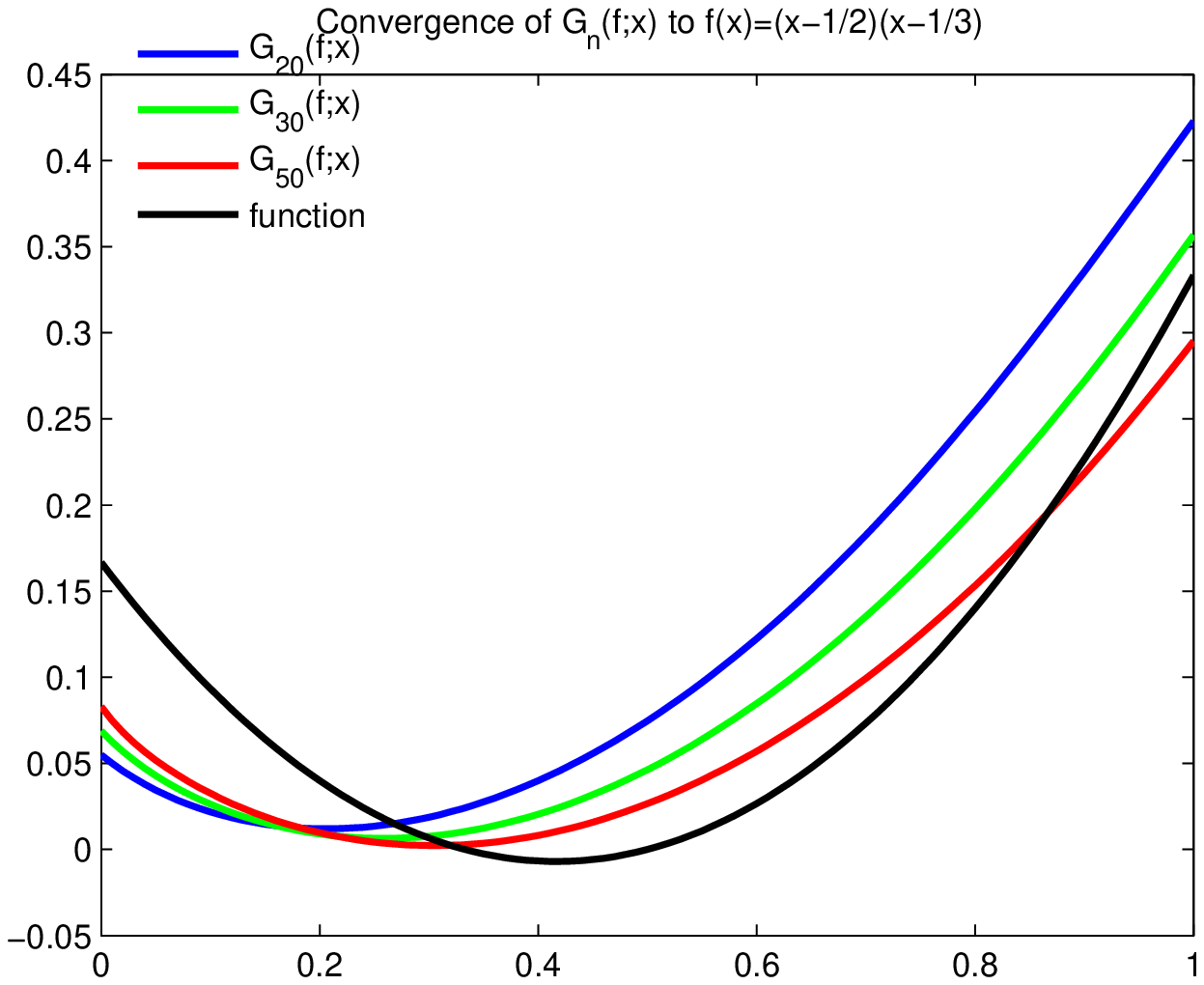, width=8cm}
\end{center}
\caption{Convergence of the operators \eqref{Gxk} to $f(x)=\left(x-\frac{1}{2}\right)\left(x-\frac{1}{3}\right)$}
\end{figure}

\begin{figure}[htb]

\begin{center}
\epsfig{file=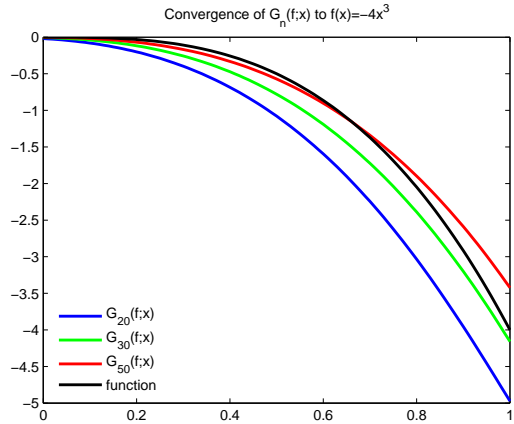, width=8cm}
\end{center}
\caption{Convergence of the operators \eqref{Gxk} to $f(x)=-4x^{3}$}
\end{figure}
It is easily seen that as the value of $n$ increases, the graphs of the operators
$G_{n}(f;x)$ converge to the graph of the function $f$.

We now compute the error estimation by using modulus of continuity for operators \eqref{Gxk}  to functions \eqref{f1} and \eqref{f2} with the help of Matlab.
Making use of the expression \eqref{G1L2} in Lemma \ref{Gav}, we find that
 \begin{equation}
|G_{n}(f;x)-f(x)|\leq||f'|| \sqrt{\frac{x}{n}+\frac{6}{n^{2}}}+\frac{1}{2}||f''||~\left( \frac{x}{n}+\frac{6}{n^{2}}\right).\label{In1}
\end{equation}
From inequality \eqref{In1}, the error bounds for the approximation by the operators \eqref{Gxk} to functions \eqref{f1} and \eqref{f2} are obtained using Matlab. These error bounds are given in Tables 1 and 2.

\begin{center}
{\bf Table 1: Error bounds by operators \eqref{Gxk} to $\left(x-\frac{1}{2}\right)\left(x-\frac{1}{3}\right)$}\\
\vspace{0.2cm}
\begin{tabular}{|c|c|c|c|}
\hline
$n$&error bound at $x=0.2$ & error bound at $x=0.5$ &error bound at $x=0.8$\\
\hline
$20$&0.0885&0.0608&0.2148\\
\hline
$30$&0.0600&0.0405&0.1600\\
\hline
$50$&0.0391&0.0260&0.1144\\
\hline
\end{tabular}
\end{center}

\vspace{0.5cm}
\begin{center}
{\bf Table 2: Error bounds by operators \eqref{Gxk} to $-4x^{3}$}\\
\vspace{0.2cm}
\begin{tabular}{|c|c|c|c|}
\hline
$n$&error bound at $x=0.2$ & error bound at $x=0.5$ &error bound at $x=0.8$\\
\hline
$20$&0.1169&0.7025&1.9651\\
\hline
$30$&0.0748&0.5066&1.4795\\
\hline
$50$&0.0462&0.3535&1.0728\\
\hline
\end{tabular}
\end{center}

\vspace{0.5cm}
\begin{eg}
{\bf The case $\mathbf{\{(k_1+k_2)(x)^{k_1+k_2-1}\}_{k_1,k_2=1}^{\infty}}$}
\end{eg}

The multiple Sheffer polynomials $S_{k_1,k_2}(x)$ reduce to the polynomials $\{(k_1+k_2)(x)^{k_1+k_2-1}\}_{k_1,k_2=1}^{\infty}$ if we consider $A(t_1,t_2)=t_1+t_2$ and $H(t_1,t_2)=t_1+t_2$,  defined by the generating function:
\begin{equation}\label{kxk}
(t_1+t_2)e^{x(t_1+t_2)}=\sum\limits_{k_1=1}^{\infty}\sum\limits_{k_2=1}^{\infty}(k_1+k_2)(x)^{k_1+k_2-1}\frac{t_{1}^{k_1}t_{2}^{k_2}}{k_{1}!~k_{2}!}.
\end{equation}

It is clear that $(k_1+k_2)(x)^{k_1+k_2-1}\ge 0$ for $x\in (0,\infty]$, $A(1,1)=2\neq 0$ and $H_{t_1}(1,1)=H_{t_2}(1,1)=1$. Hence, the restrictions \eqref{restrictions} are satisfied.
The operators \eqref{MultipleShef-Operator} involving $\{(k_1+k_2)(x)^{k_1+k_2-1}\}_{k_1,k_2=1}^{\infty}$ are given as:
\begin{equation}\label{Gkxk}
G_{n}(f;x)=\frac{e^{-nx}}{2}
\sum\limits_{k_1=1}^{\infty}\sum\limits_{k_2=1}^{\infty}
\frac{\left(\frac{nx}{2}\right)^{k_1+k_2-1}}{k_{1}!~k_{2}!}
f\left(\frac{k_{1}+k_2}{n}\right).
\end{equation}

\begin{lem}\label{G2L1}
The operators $G_{n}(f;x)$ given by \eqref{Gkxk} and $x\in [0,\infty)$, satisfy
\begin{eqnarray}
  G_{n}(e_0;x) &=& 1 \\
 G_{n}(e_1;x) &=& x+\frac{1}{n} \\
 G_{n}(e_2;x) &=& x^{2}+\frac{3x}{n}+\frac{1}{n^{2}}.
\end{eqnarray}
\end{lem}

\begin{lem}
For the operators $G_{n}(f;x)$ given by \eqref{Gkxk} and $x\in [0,\infty)$, the following hold:
\begin{eqnarray}
 G_{n}((e_1-1);x) &=& \frac{1}{n},\\
 G_{n}((e_1-x)^{2};x) &=&  \frac{x}{n}+\frac{1}{n^{2}}.\label{G2L2}
\end{eqnarray}
\end{lem}

In Figs. 3 and 4 the convergence of the operators \eqref{Gkxk} to the functions \eqref{f1} and \eqref{f2}, respectively, is illustrated.

\newpage
\begin{figure}[htb]

\begin{center}
\epsfig{file=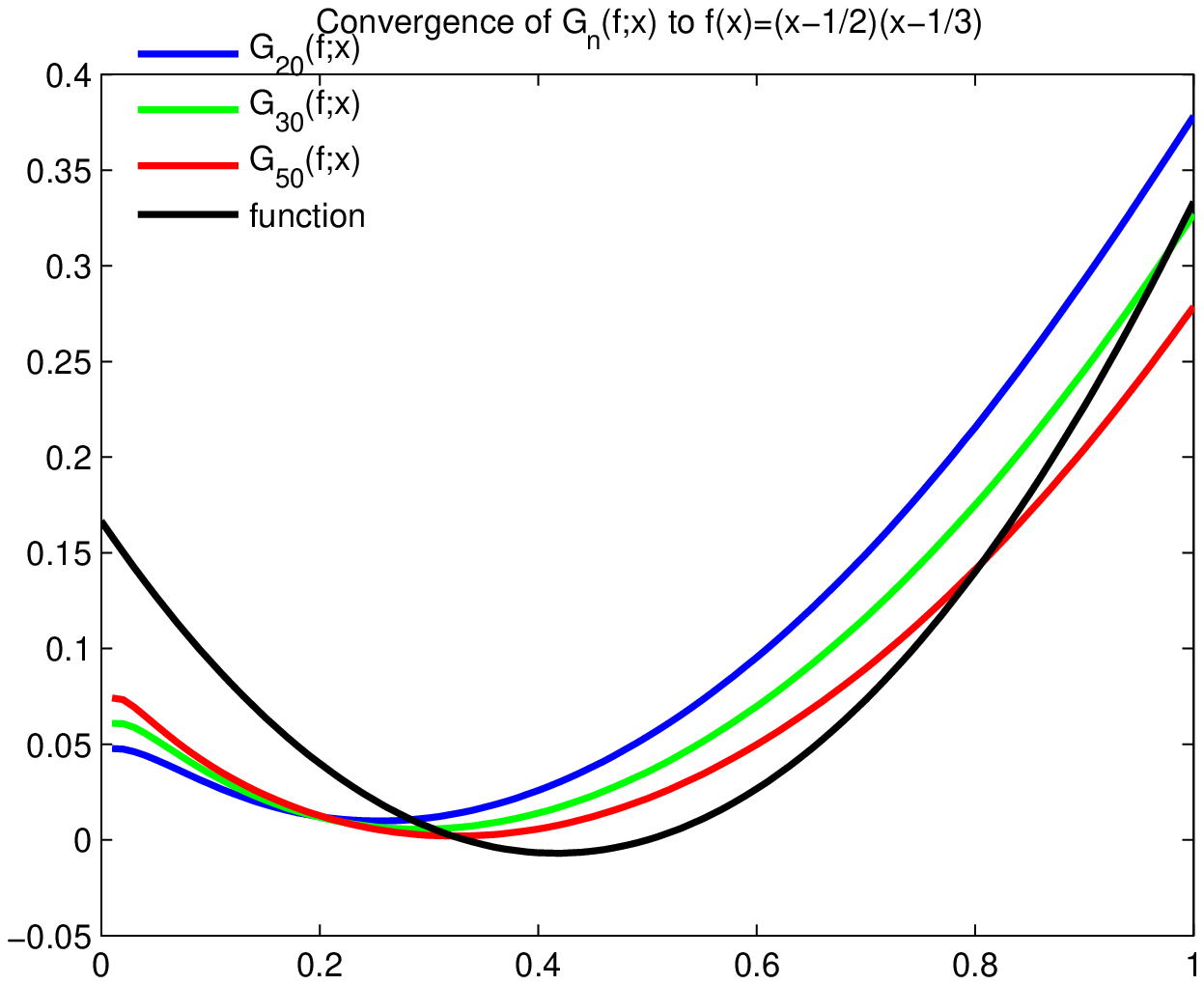, width=8cm}
\end{center}
\caption{Convergence of the operators \eqref{Gkxk} to $f(x)=\left(x-\frac{1}{2}\right)\left(x-\frac{1}{3}\right)$}
\end{figure}

\begin{figure}[htb]

\begin{center}
\epsfig{file=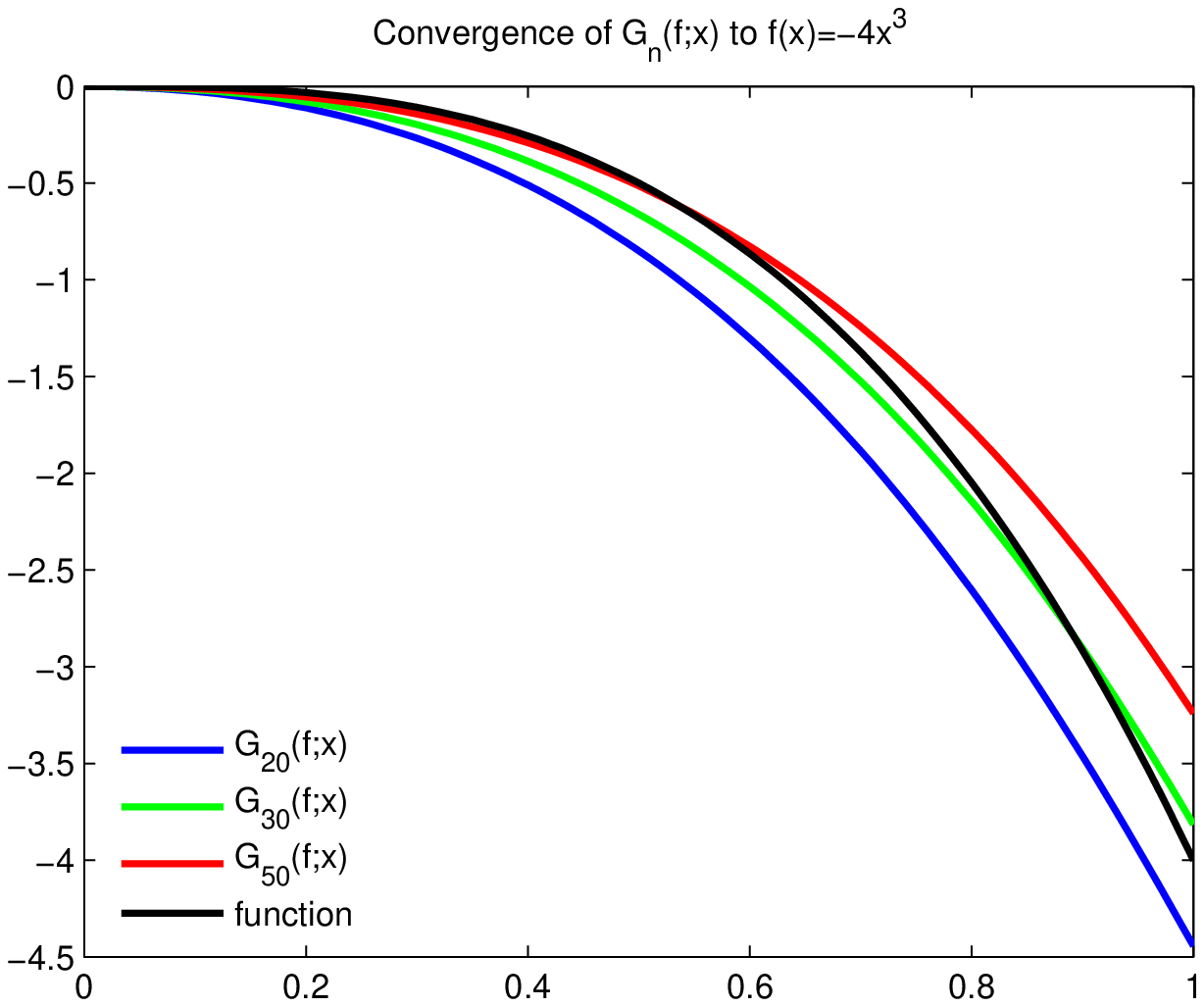, width=8cm}
\end{center}
\caption{Convergence of the operators \eqref{Gkxk} to $f(x)=-4x^{3}$}
\end{figure}

The following holds for \eqref{Gkxk}:
 \begin{equation}
|G_{n}(f;x)-f(x)|\leq||f'|| \sqrt{\frac{x}{n}+\frac{1}{n^{2}}}+\frac{1}{2}||f''||~\left( \frac{x}{n}+\frac{1}{n^{2}}\right)\label{In2}
\end{equation}

From inequality \eqref{In2}, the error bounds for the approximation by operators \eqref{Gkxk} are given in Tables 3 and 4. Inspection of Tables 1 -- 4 shows that the absolute error decreases with increasing $n$.

\begin{center}
{\bf Table 3: Error bounds by operators \eqref{Gkxk} to $\left(x-\frac{1}{2}\right)\left(x-\frac{1}{3}\right)$}\\
\vspace{0.2cm}
\begin{tabular}{|c|c|c|c|}
\hline
$n$&error bound at $x=0.2$ & error bound at $x=0.5$ &error bound at $x=0.8$\\
\hline
$20$&0.0559&0.0426&0.1806\\
\hline
$30$&0.0427&0.0317&0.1422\\
\hline
$50$&0.0311&0.0224&0.1066\\
\hline
\end{tabular}
\end{center}

\vspace{0.5cm}
\begin{center}
{\bf Table 4: Error bounds by operators \eqref{Gkxk} to $-4x^{3}$}\\
\vspace{0.2cm}
\begin{tabular}{|c|c|c|c|}
\hline
$n$&error bound at $x=0.2$ & error bound at $x=0.5$ &error bound at $x=0.8$\\
\hline
$20$&0.0647&0.5250&1.6273\\
\hline
$30$&0.0483&0.4150&1.3040\\
\hline
$50$&0.0348&0.3133&0.9954\\
\hline
\end{tabular}
\end{center}

\vspace{0.2cm}
\noindent
{\bf Acknowledgement.}~This work has been sponsored by a Dr. D. S. Kothari Post Doctoral Fellowship (Award letter No. F.4-2/2006(BSR)/MA/17-18/0025) awarded to {\bf Dr. Mahvish Ali} by the University Grants Commission, Government of India, New Delhi.

\end{document}